\documentclass[reqno]{amsart}

\usepackage{amssymb}
\usepackage{graphicx}
\usepackage{amscd}
\usepackage[hidelinks]{hyperref}
\usepackage{color}
\usepackage{float}
\usepackage{graphics,amsmath,amssymb}
\usepackage{amsthm}
\usepackage{amsfonts}
\usepackage{latexsym}
\usepackage{epsf}
\usepackage{enumerate}
\usepackage{xifthen}
\usepackage{mathrsfs}
\usepackage{dsfont}
\usepackage{makecell}
\usepackage[FIGTOPCAP]{subfigure}
\usepackage{amsmath}
\allowdisplaybreaks[4]
\usepackage{listings}
\usepackage{etoolbox}
\usepackage{fancyhdr}
\usepackage{pdflscape}
\usepackage[title,toc,titletoc]{appendix}
\usepackage{enumitem}
\usepackage[noadjust]{cite}

 \newtheoremstyle{mytheorem}
 {3pt}
 {3pt}
 {\slshape}
 {}
 {\bfseries}
 {.}
 { }
 {}

\numberwithin{equation}{section}

\theoremstyle{theorem}
\newtheorem{theorem}{Theorem}[section]

\newtheorem{lemma}[theorem]{Lemma}
\newtheorem{proposition}[theorem]{Proposition}

\theoremstyle{definition}

\newcommand{\Keywords}[1]{\ifthenelse{\isempty{#1}}{}{\smallskip \smallskip \noindent \textbf{Keywords}. #1}}
\newcommand{\MSC}[2][2010]{\ifthenelse{\isempty{#2}}{}{\smallskip \smallskip \noindent \textbf{#1MSC}. #2}}
\newcommand{\abstractnote}[1]{\ifthenelse{\isempty{#1}}{}{\smallskip \smallskip \noindent \textsuperscript{\dag}#1}}

\makeatletter
\def\specialsection{\@startsection{section}{1}%
  \z@{\linespacing\@plus\linespacing}{.5\linespacing}%
  {\normalfont}}
\def\section{\@startsection{section}{1}%
  \z@{.7\linespacing\@plus\linespacing}{.5\linespacing}%
  {\normalfont\scshape}}
\patchcmd{\@settitle}{\uppercasenonmath\@title}{\Large\boldmath}{}{}
\patchcmd{\@settitle}{\begin{center}}{\begin{flushleft}}{}{}
\patchcmd{\@settitle}{\end{center}}{\end{flushleft}}{}{}
\patchcmd{\@setauthors}{\MakeUppercase}{\normalsize}{}{}
\patchcmd{\@setauthors}{\centering}{\raggedright}{}{}
\patchcmd{\section}{\scshape}{\large\bfseries\boldmath}{}{}
\patchcmd{\subsection}{\bfseries}{\bfseries\boldmath}{}{}
\renewcommand{\@secnumfont}{\bfseries}
\patchcmd{\@startsection}{\@afterindenttrue}{\@afterindentfalse}{}{}
\patchcmd{\abstract}{\leftmargin3pc}{\leftmargin1pc}{}{}

\def\maketitle{\par
  \@topnum\z@ 
  \@setcopyright
  \thispagestyle{empty}
  \ifx\@empty\shortauthors \let\shortauthors\shorttitle
  \else \andify\shortauthors
  \fi
  \@maketitle@hook
  \begingroup
  \@maketitle
  \toks@\@xp{\shortauthors}\@temptokena\@xp{\shorttitle}%
  \toks4{\def\\{ \ignorespaces}}
  \edef\@tempa{%
    \@nx\markboth{\the\toks4
      \@nx\MakeUppercase{\the\toks@}}{\the\@temptokena}}%
  \@tempa
  \endgroup
  \c@footnote\z@
  \@cleartopmattertags
}
\makeatother

\title{Note on sums involving the Euler function}

\author[S. Chern]{Shane Chern}
\address{Department of Mathematics, Penn State University, University Park, PA 16802, USA}
\email{shanechern@psu.edu}

\date{}

\begin{document}

{\footnotesize\noindent To appear in \textit{Bull. Aust. Math. Soc.}}

\bigskip \bigskip

\maketitle

\begin{abstract}

In this note, we provide refined estimates of the following sums involving the Euler totient function:
$$\sum_{n\le x} \phi\left(\left[\frac{x}{n}\right]\right) \qquad \text{and} \qquad \sum_{n\le x} \frac{\phi([x/n])}{[x/n]}$$
where $[x]$ denotes the integral part of real $x$. The above summations were recently considered by Bordell\`es et al.~and Wu.


\MSC{11A25, 11L07.}
\end{abstract}

\section{Introduction}

Let $[x]$ denote the integral part of a real number $x$. In a recent paper of Bordell\`es, Dai, Heyman, Pan and Shparlinski \cite{BDHPS2018}, the asymptotic behavior of the following function was studied:
$$S_{f}:=\sum_{n\le x} f\left(\left[\frac{x}{n}\right]\right).$$
In particular, if $f(n)$ is set to be $\phi(n)$ and $\phi(n)/n$ where $\phi(n)$ is the Euler totient function, Bordell\`es et al.~obtained the following estimates.

\smallskip

\begin{equation}\label{eq:phi-n}
\sum_{n\le x} \frac{\phi([x/n])}{[x/n]}=x\sum_{n\ge 1}\frac{\phi(n)}{n^2(n+1)}+O(x^{1/2})
\end{equation}
and
\begin{align}
&\left(\frac{2629}{4009}\cdot \frac{1}{\zeta(2)}+o(1)\right)x\log x\nonumber\\
&\qquad\qquad\qquad\le\sum_{n\le x} \phi\left(\left[\frac{x}{n}\right]\right)\le \left(\frac{2629}{4009}\cdot \frac{1}{\zeta(2)}+\frac{1380}{4009}+o(1)\right)x\log x\label{eq:phi}
\end{align}
for $x\to\infty$.

\medskip

Subsequently, Wu respectively improved in \cite{Wu2018a} the upper and lower bounds in \eqref{eq:phi} and in \cite{Wu2018b} the error term in \eqref{eq:phi-n}. More precisely, Wu showed that the error term in \eqref{eq:phi-n} can be sharpened to $O(x^{1/3}\log x)$, while the bounds in \eqref{eq:phi} can be refined as
\begin{equation}\label{eq:bound-wu}
\frac{2}{3}\cdot \frac{1}{\zeta(2)}x\log x+O(x)\le \sum_{n\le x} \phi\left(\left[\frac{x}{n}\right]\right) \le \left(\frac{2}{3}\cdot \frac{1}{\zeta(2)}+\frac{1}{3}\right)x\log x+O(x).
\end{equation}

\medskip

To bound $\displaystyle \sum_{n\le x} \phi\left(\left[\frac{x}{n}\right]\right)$, the main idea in Bordell\`es et al.~\cite{BDHPS2018} and Wu \cite{Wu2018a} relies on an estimate of the following summation
$$\mathfrak{S}_{\delta}(x,N):=\sum_{N<n\le 2N}\phi(n)\psi\left(\frac{x}{n+\delta}\right)$$
for $x\ge 2$ and $1\le N\le x$ where $\psi(x)=x-[x]-\frac{1}{2}$ and $\delta\in\{0,1\}$. Such an estimate is built on Vaaler's expansion formula of $\psi(x)$ (cf.~\cite{Vaa1985} or Theorem 6.1 in \cite{Bor2012}) and the theory of exponential pairs (cf.~Section 6.6.3 in \cite{Bor2012}). Further, as Wu has shown in \cite{Wu2018b}, the estimate of a similar summation
$$\mathfrak{S}_{\delta}^*(x,N):=\sum_{N<n\le 2N}\frac{\phi(n)}{n}\psi\left(\frac{x}{n+\delta}\right)$$
will be useful to deduce the error term in \eqref{eq:phi-n}.

\medskip

We observe that, with the aid of an elaborate result due to Huxley (cf.~\cite{Hux2003} or Theorem 6.40 in \cite{Bor2012}), the estimate of $\mathfrak{S}_{\delta}^*(x,D)$ in \cite{Wu2018b} can be further sharpened. In fact, Huxley's result is strong enough in the sense that the best known error term up to now for the Dirichlet divisor problem can be deduced from it.

\smallskip

In this note, we shall prove the following results.

\smallskip

\begin{theorem}\label{th:1}
We have, as $x\to\infty$,
\begin{equation}\label{eq:phi-n-1}
\sum_{n\le x} \frac{\phi([x/n])}{[x/n]}=x\sum_{n\ge 1}\frac{\phi(n)}{n^2(n+1)}+O(x^{131/416}\left(\log x\right)^{26947/8320}).
\end{equation}
\end{theorem}

\smallskip

\begin{theorem}\label{th:2}
We have, as $x\to\infty$,
\begin{align}
&\frac{285}{416}\cdot \frac{1}{\zeta(2)}x \log x+O(x \log\log x)\nonumber\\
&\qquad\qquad\qquad\le\sum_{n\le x} \phi\left(\left[\frac{x}{n}\right]\right)\le \left(\frac{285}{416}\cdot \frac{1}{\zeta(2)}+\frac{131}{416}\right)x \log x+O(x \log\log x).\label{eq:phi-1}
\end{align}
\end{theorem}

\medskip

We have two remarks to make.

\begin{enumerate}[label=\arabic*.]
\item Let $\tau(n)$ denote the number of divisors of $n$. It is known that the main term of $\displaystyle \sum_{n\le x}\tau(n)$ is $x(\log x + 2\gamma -1)$ where $\gamma$ is the Euler constant. The error term, denoted by $\Delta(x)$, can be trivially bounded to be $O(x^{1/2})$. Hardy \cite{Har1916} also showed that $\Delta(x)$ cannot be $o(x^{1/4})$. The best known bound up to now for $\Delta(x)$ is $O(x^{131/416}\left(\log x\right)^{26947/8320})$ which is due to Huxley as we have mentioned above. We can see that the error term in \eqref{eq:phi-n-1} can also reach this size.

\item Numerically, we have
\begin{align*}
\frac{285}{416}\cdot \frac{1}{\zeta(2)}\approx 0.41649 \qquad &\text{and} \qquad \frac{285}{416}\cdot \frac{1}{\zeta(2)}+\frac{131}{416}\approx 0.73139.\\
\intertext{This slightly improves the bounds of Wu \cite{Wu2018a} in \eqref{eq:bound-wu}:}
\frac{2}{3}\cdot \frac{1}{\zeta(2)}\approx 0.40528 \qquad &\text{and} \qquad \frac{2}{3}\cdot \frac{1}{\zeta(2)}+\frac{1}{3}\approx 0.73861.
\end{align*}
\end{enumerate}

\section{An auxiliary estimate}

Let $\delta\in\{0,1\}$. We will focus on the following auxiliary function defined in the Introduction:
$$\mathfrak{S}_{\delta}^*(x,N):=\sum_{N<n\le 2N}\frac{\phi(n)}{n}\psi\left(\frac{x}{n+\delta}\right).$$

One has
\begin{align}
\sum_{N<n\le 2N}\frac{\phi(n)}{n}\psi\left(\frac{x}{n+\delta}\right)&=\sum_{N<n\le 2N}\frac{1}{n}\psi\left(\frac{x}{n+\delta}\right) \sum_{\substack{k,\ell\\k\ell=n}}\mu(k)\ell\nonumber\\
&=\sum_{k\le 2N}\frac{\mu(k)}{k}\sum_{N/k<\ell\le 2N/k}\psi\left(\frac{x}{k\ell+\delta}\right).\label{eq:sum_phi_n_psi}
\end{align}
Now we will apply the following result due to Huxley \cite{Hux2003}.

\begin{lemma}[Huxley, cf.~Theorem 6.40 in \cite{Bor2012}]\label{le:main}
Let $r\ge 5$, $M\ge 1$ be integers and $f\in \mathcal{C}^{r}[M,2M]$ such that there exist real numbers $T\ge 1$ and $1\le c_0\le \cdots\le c_r$ such that, for all $x\in[M,2M]$ and all $j\in\{0,\ldots,r\}$, we have
$$\frac{T}{M^j}\le |f^{(j)}(x)| \le c_j \frac{T}{M^j}.$$
Then we have
$$\sum_{M< n\le 2M}\psi\big(f(n)\big)\ll (MT)^{131/416}(\log MT)^{18627/8320}.$$
\end{lemma}

\medskip

Under the setting of Lemma \ref{le:main}, let us put $\displaystyle f(z)=\frac{x}{kz+\delta}$. It can be easily computed that for $j\ge 1$
$$f^{(j)}(z)=(-1)^j\frac{j!k^jx}{(kz+\delta)^{j+1}}.$$
Notice that we have, trivially, $\displaystyle \left[\frac{N}{k}\right]\asymp \frac{N}{k}$ when $k< N$. It can be shown with almost no effort that $T$ can be chosen to be $\displaystyle \asymp \frac{x}{N}$. In fact, $\displaystyle T=C\frac{x}{N}$ is admissible where $\displaystyle C=\frac{6!}{3\cdot 6^6}$. Now we assume that $N\le Cx$.

\medskip

For $k< N$, we notice that
$$\sum_{N/k<\ell\le 2N/k}\psi\left(\frac{x}{k\ell+\delta}\right) = \sum_{[N/k]<\ell\le 2[N/k]}\psi\left(\frac{x}{k\ell+\delta}\right)+O(1).$$
It follows that, for $k< N\le Cx$,
$$\sum_{N/k<\ell\le 2N/k}\psi\left(\frac{x}{k\ell+\delta}\right)\ll \left(\frac{x}{k}\right)^{131/416}\left(\log \frac{Cx}{k}\right)^{18627/8320}.$$

\medskip

Further, for $N\le k\le 2N$,
$$\sum_{N/k<\ell\le 2N/k}\psi\left(\frac{x}{k\ell+\delta}\right)\ll 1.$$

\medskip

Hence, by \eqref{eq:sum_phi_n_psi}, we conclude that
\begin{align*}
\sum_{N<n\le 2N}\frac{\phi(n)}{n}\psi\left(\frac{x}{n+\delta}\right)&=\sum_{k\le 2N}\frac{\mu(k)}{k}\sum_{N/k<\ell\le 2N/k}\psi\left(\frac{x}{k\ell+\delta}\right)\\
&\ll \sum_{k< N}\frac{1}{k}\left(\frac{x}{k}\right)^{131/416}\left(\log x\right)^{18627/8320}\\
&\ll x^{131/416}\left(\log x\right)^{18627/8320}.
\end{align*}

\medskip

To summarize, we have

\begin{proposition}\label{prop:main}
Let $\delta\in\{0,1\}$. Then
\begin{equation}
\sum_{N<n\le 2N}\frac{\phi(n)}{n}\psi\left(\frac{x}{n+\delta}\right)\ll x^{131/416}\left(\log x\right)^{18627/8320}
\end{equation}
uniformly for $\displaystyle 1\le N\le \frac{6!}{3\cdot 6^6}x$.
\end{proposition}

\section{Proof of Theorem \ref{th:1}}

Again, let $\displaystyle C=\frac{6!}{3\cdot 6^6}$. Following the argument in \cite{Wu2018b}, we have
\begin{equation}
\sum_{n\le x} \frac{\phi([x/n])}{[x/n]}=\sum_{1/C<n\le x} \frac{\phi([x/n])}{[x/n]}+O(1).
\end{equation}

\medskip

Notice that if $d=[x/n]$, then $x/(d+1)<n\le x/d$. Hence,
\begin{align*}
\sum_{1/C<n\le x} \frac{\phi([x/n])}{[x/n]}&=\sum_{d\le Cx} \frac{\phi(d)}{d} \sum_{x/(d+1)<n\le x/d}1\\
&=\sum_{d\le Cx} \frac{\phi(d)}{d} \Bigg(\frac{x}{d(d+1)}+\psi\left(\frac{x}{d+1}\right)-\psi\left(\frac{x}{d}\right)\Bigg)\\
&=x\sum_{d\ge 1}\frac{\phi(d)}{d^2(d+1)}+O(1)+\sum_{d\le Cx}\frac{\phi(d)}{d}\Bigg(\psi\left(\frac{x}{d+1}\right)-\psi\left(\frac{x}{d}\right)\Bigg).
\end{align*}
Using a dyadical split together with Proposition \ref{prop:main}, we see that for $\delta\in\{0,1\}$,
$$\sum_{d\le Cx}\frac{\phi(d)}{d}\psi\left(\frac{x}{d+\delta}\right)\ll x^{131/416}\left(\log x\right)^{26947/8320}.$$

\medskip

We therefore arrive at Theorem \ref{th:1}.

\section{Proof of Theorem \ref{th:2}}

We first split the sum $\displaystyle \sum_{n\le x} \phi\left(\left[\frac{x}{n}\right]\right)$ into two parts:
\begin{equation}
\sum_{n\le x} \phi\left(\left[\frac{x}{n}\right]\right)=\sum_{n\le D}\phi\left(\left[\frac{x}{n}\right]\right)+\sum_{D<n\le x}\phi\left(\left[\frac{x}{n}\right]\right),
\end{equation}
where $D\ge 1/C$ is to be determined later.

\medskip

Using a similar argument to that in the previous section, we have
\begin{align}
\sum_{D<n\le x} \phi\left(\left[\frac{x}{n}\right]\right)&=\sum_{d\le x/D} \phi(d) \sum_{x/(d+1)<n\le x/d}1\nonumber\\
&=x\sum_{d\le x/D}\frac{\phi(d)}{d(d+1)}+\sum_{d\le x/D} \phi(d)\Bigg(\psi\left(\frac{x}{d+1}\right)-\psi\left(\frac{x}{d}\right)\Bigg)\nonumber\\
&=\frac{1}{\zeta(2)} x \log \frac{x}{D}+O(x)+\sum_{d\le x/D} \phi(d)\Bigg(\psi\left(\frac{x}{d+1}\right)-\psi\left(\frac{x}{d}\right)\Bigg).\label{eq:D-x}
\end{align}
Here in the last identity we use the following standard result (cf.~Exercise 3.6 in \cite{Apo1976})
$$\sum_{n\le x}\frac{\phi(n)}{n^2}=\frac{1}{\zeta(2)}\log x+O(1)$$
so that
$$\sum_{n\le x}\frac{\phi(n)}{n(n+1)}=\sum_{n\le x}\phi(n)\Bigg(\frac{1}{n^2}+O\left(\frac{1}{n^3}\right)\Bigg)=\frac{1}{\zeta(2)}\log x+O(1).$$

\medskip

Applying Abel's summation formula to the last part in \eqref{eq:D-x} yields
\begin{align}
\sum_{d\le x/D} \phi(d)\Bigg(\psi\left(\frac{x}{d+1}\right)-\psi\left(\frac{x}{d}\right)\Bigg)&=\frac{x}{D}\sum_{d\le x/D}\frac{\phi(d)}{d}\Bigg(\psi\left(\frac{x}{d+1}\right)-\psi\left(\frac{x}{d}\right)\Bigg)\nonumber\\
&\quad-\int_{1}^{x/D}\sum_{d\le t}\frac{\phi(d)}{d}\Bigg(\psi\left(\frac{x}{d+1}\right)-\psi\left(\frac{x}{d}\right)\Bigg) dt.\label{eq:abel-sum}
\end{align}

\medskip

Notice that for $t\in [1,x/D]$ and $\delta\in\{0,1\}$, by a dyadical split, it follows from Proposition \ref{prop:main} that
\begin{align*}
\sum_{d\le t}\frac{\phi(d)}{d}\psi\left(\frac{x}{d+\delta}\right)\ll x^{131/416}\left(\log x\right)^{26947/8320}.
\end{align*}
It turns out that by \eqref{eq:abel-sum}
\begin{align}
\sum_{d\le x/D} \phi(d)\Bigg(\psi\left(\frac{x}{d+1}\right)-\psi\left(\frac{x}{d}\right)\Bigg)\ll \frac{x}{D}x^{131/416}\left(\log x\right)^{26947/8320}.\label{eq:D-x-error}
\end{align}

\medskip

Let us choose
$$D=x^{131/416}\left(\log x\right)^{26947/8320}.$$

\medskip

It follows from \eqref{eq:D-x} and \eqref{eq:D-x-error} that
\begin{align}
\sum_{D<n\le x}\phi\left(\left[\frac{x}{n}\right]\right)&=\frac{1}{\zeta(2)}\left(1-\frac{131}{416}\right)x \log x+O(x \log\log x)\nonumber\\
&=\frac{285}{416}\cdot \frac{1}{\zeta(2)}x \log x+O(x \log\log x).\label{eq:2-bound1}
\end{align}
We can also trivially bound
\begin{align}
0\le \sum_{n\le D}\phi\left(\left[\frac{x}{n}\right]\right)&\le \sum_{n\le D}\frac{x}{n}\nonumber\\
&=\frac{131}{416} x\log x+O(x \log\log x).\label{eq:2-bound2}
\end{align}

\medskip

Theorem \ref{th:2} is a direct combination of \eqref{eq:2-bound1} and \eqref{eq:2-bound2}.

\subsection*{Acknowledgements}

I want to thank Jie Wu for sharing the manuscript of \cite{Wu2018b}.

\bibliographystyle{amsplain}

\end{document}